\documentclass[prc,superscriptaddress,mathtools]{revtex4}
\usepackage{natbib}
\usepackage{amsmath}
\usepackage{amssymb}
\usepackage{color}
\usepackage{graphicx}
\usepackage{bm}

\usepackage{multirow}

\begin{document}

\title{Symmetry analysis of an elastic beam with axial load}

\author{Bidisha Kundu$^*$, Ranjan Ganguli\\
Department of Aerospace Engineering, Indian Institute of Science,
Bangalore 560012, India\\ bidishakundu@iisc.ac.in, ganguli@iisc.ac.in\\
$^*$Corresponding author
}




\begin{abstract}
We construct the closed form solution of an elastic beam with axial load using Lie symmetry method. A beam with spatially varying physical properties such as mass and second moment of inertia is considered. The governing fourth order partial differential equation with variable coefficients which is not amenable to simple methods of solution, is solved using Lie symmetry. We incorporate boundary conditions and then compare with the numerical solution.
\end{abstract}

\maketitle
\section{Introduction}
An elastic beam is a three dimensional structure whose axial extension is more than any other dimension orthogonal to it. It is a fundamental model which pervades every corner of physics and engineering.
However, the integrability in finite numbers of terms or exact solution of the governing equation of this model is an open question.
The deflection is dependent not only on the mass of the beam or external applied force but also on the elastic properties of the material and geometry of the beam.
The equation of motion of the beam is a fourth-order linear partial differential equation with variable coefficients. The equation can be derived from generalised Hooke's law and force balance or by minimizing the energy of the system. Though the equation is linear in nature, due to the presence of variable coefficients, it is very arduous to get the solution analytically. In \cite{guede2001apparentlynew}, the closed form solution of vibrating axially loaded beam was studied with different boundary conditions. However, to the best of our knowledge, in this work, for the first time, Euler-Bernoulli beam with axial force is studied using the Lie-method. Here we follow the Lie-group method to get the similarity solution.
\par
The main motivation to use the Lie method is to get the solution of the mathematical model in different forms which are very easy to use. First, Sophus Lie applied this method to partial differential equations and later this method was further developed by Ovsjannikov \cite{ovsjannikov1962group} and Matschat and M\"{u}ller\cite{muller1962and}.

Lie method has been applied to various types of ordinary differential equations (ODE) and partial differential equations(PDE). A systematic approach for applying Lie method to ODE and PDE is found in Refs. [\cite{olver2012applications},\cite{bluman2013symmetries},\cite{ibragimov1995crc}, \cite{hydon2000symmetrynew}]. Bluman et al. used this method for different types of
mathematical physics
problems such as wave equation, diffusion equation etc. \cite{bluman1980remarkable,bluman1987invariance}. Torrisi et al. studied diffusion equation by equivalence transformation \cite {torrisi1996group}. Ibragimov applied this method to some real life problems \cite{ibragimov2004equivalence,ibragimov2011lie} related to tumour growth model and also to metallurgical industry mathematical models.
In \cite{gray2015calculate}, the procedure to calculate all point symmetries of linear and linearizable differential equations was studied. An algorithm for integrating systems of two second-order ordinary differential equations with four symmetries was given in \cite{gainetdinova2017integrability}. The Lie symmetry method is also employed in ordinary difference equation \cite{hydon2000symmetries}.

\par
The Lie-method is also available in recent research addressing equations in mathematical physics \cite{kang2012symmetry,singla2017invariant,hau2017optimal}.
In Ref. \cite{kang2012symmetry}, Kang and Qu have studied the relationship between Lie point symmetry and fundamental solution for systems of parabolic equations.
The invariant analysis of space-time fractional nonlinear systems of partial differential equations is performed using Lie symmetry method in Ref. \cite{singla2017invariant}.
Hau, et al. presented a unifying solution framework for the linearized compressible equations for two-dimensional linearly sheared unbounded flows \cite{hau2017optimal}.

The symmetry analysis for a physical problem is very important. In Lorentz transformation used in special relativity, Yang-Mills theory \cite{strobl2004algebroid,beisert2017yangian}, and Schr\"{o}dinger equation, the Lie symmetry analysis has been used extensively. In Ref. \cite{belmonte2007lie}, Lie symmetries and canonical transformations are applied to construct the explicit solutions of Schr\"{o}dinger equation with a spatially inhomogeneous nonlinearity from  those of the homogeneous nonlinear Schr\"{o}dinger equation.
In Ref. \cite{budanur2015reduction}, symmetry analysis is used to understand the fluid flow in a pipe or channel structure.
\par
The Lie-method has also been employed in the field of mechanics \cite{bocko2012symmetries}.
In Ref. \cite{ozkaya2002group}, a beam moving with time-dependent axial velocity is studied using equivalence transformation. Wafo discussed the Euler-Bernoulli (EB) beam from a symmetry stand-point \cite{soh2008euler}.
The general beam equation is studied by Bokhari et al. for symmetries and integrability with Lie method \cite{bokhari2010symmetries}. Bokhari et al also found the complete Lie symmetry classification of the fourth-order dynamic EB beam equation with load dependent on normal displacement. Johnpillai et al studied the EB beam equation from the Noether symmetry viewpoint \cite{johnpillai2016noether}.
\par
 We search for one-parameter group of transformation which leaves the governing PDE invariant and we get the corresponding Lie-algebras. Using the infinitesimal generator of this transformation, we solve the newly reduced differential equation. In Section $2$, we describe our problem. The mathematical theories for the Lie-symmetry approach, the procedure to get the invariant solution and application of the Lie method to our problem are described in Section $3$. The closed form solutions for different cases and boundary condition imposition are given in Section $4$ and Section $5$, respectively.

\section{Formulation}
We fix the coordinate axes $X$, $Y, Z$ along the length, breadth and height of the beam, respectively. We consider that the beam is slender and
of length $l$. Here, $u(x,t)$ and $M(x,t)$ are the out-plane bending displacement along Z and bending moment, respectively,
at the point $x$ and the instant $t$ where $x\in [0,l]$ and $t\in \mathbb{R}^{+}$.
We assume the stiffness functions $EI:[0,l]\rightarrow \mathbb{R}^{+}$ and $m:[0,l]\rightarrow \mathbb{R}^{+}$ are continuously differentiable functions. Here $\mathbb{R}^{+}$ is the set of all positive real numbers.
\par
\begin{figure}[h]
\centering
\includegraphics[scale=0.4]{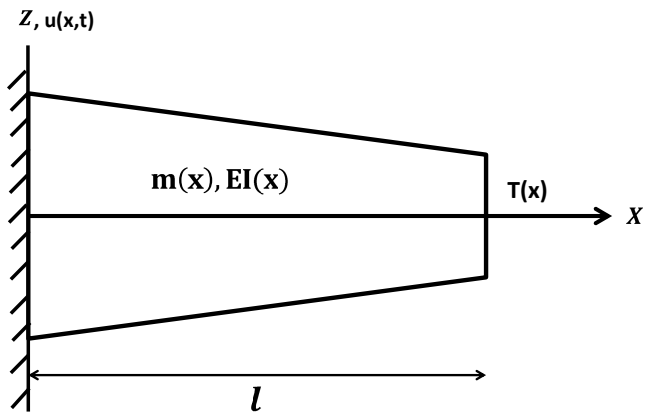}
\caption{Schematic diagram of an axially loaded beam (under tension).}
\label{beam}
\end{figure}

The equation of motion of the Euler-Bernoulli beam with axial force is given by
\begin{eqnarray}
\label{first equation}
\frac{\partial{^2}}{\partial x{^2}}\left(EI(x)\frac{\partial{^2}u}{\partial x{^2}}\right)+
m(x)\frac{\partial{^2}u}{\partial t{^2}}-
\frac{\partial}{\partial x}\left(T(x)\frac{\partial u}{\partial x}\right)=0.
\end{eqnarray}
 For rotating beam $T(x)=\int_{x}^{l}m(x)\Omega^{2}x dx$ where $\Omega$ is the rotating speed. Another important practical case is the gravity loaded beam. For this problem $T(x)=\int_{x}^{l}m(x)g dx$ where $g$ is the gravitational force acting on the system. In the case of a stiff-string or piano string $T$ is constant \cite{gunda2008stiff}.
 This formulation is valid for all tensile loads and for compressive loads less than the critical buckling load.
In general, for this type of beam, cantilever boundary condition at the left end and free at the right end is considered.
Here, a schematic diagram of an axially loaded beam (under tension), is
shown in Figure (\ref{beam}).
\vspace*{-7pt}
\vspace*{-5pt}
\section{The Lie symmetry method}
Eq. (\ref{first equation}) is a linear PDE of order four in $x$ and order two in $t$.
The Lie symmetry method will help us to extract the closed form solution via the symmetry analysis. We define the original Eq. (\ref{first equation}) in geometric form.
\par
\par
Consider $X=\{(x,t)|x,t\in \mathbb{R}\}$  as the two-dimensional manifold of the independent variables and $U\subset \mathbb{R}$ as the one dimensional manifold of dependent variable.
Equation (\ref{first equation}) is defined on the manifold $\mathbb{M}\subset X\times U$ of dimension three.
In order to establish a concrete geometric structure of the differential equation, we have to construct a manifold $\mathfrak{J}$ of dimension $p$ which includes the independent variables, dependent variables,
and all possible partial derivatives of the dependent variable with respect to independent variables up to order of the given differential equation. Here the dimension $p$ depends on the dimension of $X, U$ as well as on the order of the differential equation given. We ``prolong'' the original manifold $\mathbb{M}\subset X\times U$ to $\mathfrak{J}$ which also has manifold structure.
\par
The main motivation to study a differential equation with Lie symmetry method is to find a favourable coordinate transformation which can produce another equivalent form of the given equation. In the new coordinate systems,
the new equation may have closed form solution. To get the favourable coordinate transformations, we should take care of the derivatives of the dependent variables which is related to the vector fields of $\mathbb{M}$. We also want to measure the contribution of the vector fields related to the given differential equation in the prolonged space $\mathfrak{J}$ which is called the prolongation of the vector field.
\par
As this is a PDE of order four, we need the prolongation of order four, i.e.
$pr^{(4)}{\bf v}$.
For a PDE with one dependent variable and two independent variables, the one-parameter Lie group of transformations is
\begin{eqnarray}
\label{one-parameter}
 x^{*}=X(x,t,u;\epsilon )=  x+\epsilon  {\xi}(x,t,u) + O (\epsilon^{2}) \nonumber \\
  t^{*}=T(x,t,u;\epsilon )=  t+\epsilon \tau(x,t,u) + O (\epsilon^{2})\nonumber \\
  u^{*}=U(x,t,u;\epsilon )=  u+\epsilon  {\eta}(x,t,u) + O (\epsilon^{2})
 \end{eqnarray}
 The prolongation of order four $pr^{(4)}{\bf v}$ is given by
 \begin{eqnarray}
pr^{(4)}{\bf v}={\xi}\frac{\partial}{\partial x}+\tau \frac{\partial}{\partial t}+{\eta} \frac{\partial}{\partial u}
+{\eta^x} \frac{\partial}{\partial u_{x}}+{\eta^t} \frac{\partial}{\partial u_{t}}+{\eta^{xx}} \frac{\partial}{\partial u_{xx}}+...
\nonumber \\
{\eta^{xxt}} \frac{\partial}{\partial u_{xxt}}+...+{\eta^{ttt}} \frac{\partial}{\partial u_{ttt}}
+{\eta^{xxxx}} \frac{\partial}{\partial u_{xxxx}}+{\eta^{xxxt}} \frac{\partial}{\partial u_{xxxt}}+...+{\eta^{tttt}}\frac{\partial}{\partial u_{tttt}}
\end{eqnarray}
\par
 From the Fundamental theorem [Theorem $2.31$, \cite{olver2012applications}] of transformations (\ref{one-parameter}) admitted by the PDE (\ref{first equation}), applying $pr^{(4)}{\bf v}$ on Eq.(\ref{first equation}) yields
\begin{gather}
 \xi (x,t,u(x,t)) \left( EI^{(3)}(x) u_{xx}(x,t)+2 EI''(x)
   u_{xxx}(x,t)\right.\nonumber\\
    \left.
   +EI'(x) u_{xxxx}(x,t)+
   m'(x) u_{tt}(x,t)
  - T''(x)u_{x}(x,t)\right.\nonumber\\
    \left.
 - T'(x) u_{xx}(x,t)  \right) +
   \eta^{ xx} EI''(x)+2
   \eta ^{xxx} EI'(x)\nonumber\\+\eta^{xxxx} EI(x)+
   \eta ^{tt}m(x)
   -\eta^{x} T'(x)- \eta ^{xx} T(x)=0
 \label{equaprolong}
\end{gather}
Eq. (\ref{equaprolong}) should be satisfied to yield the admissible transformations.
The value of $u_{xxxx}(x,t)$ from Eq. (\ref{first equation}) is substituted into the above Eq. (\ref{equaprolong}).
 Now substituting the values of $\eta^{x}, \eta^{xx}, \eta^{tt}, \eta^{xxx}, \eta^{xxxx}$ from
  Theorem 32.3.5 \cite{hassani2013mathematical},
Eq.
(\ref{equaprolong}) is simplified
 and a polynomial equation is formed in all possible
derivatives $u^{i,j}$ of $u$ upto order four with respect to $x$, $t$ where
\begin{eqnarray}
 u^{i,j}=\frac{\partial^{i+j} u}{\partial x^{i}\partial t^{j}}
\end{eqnarray}
This polynomial equation should be satisfied for arbitrary $x$, $t$ and $u^{i,j}$ which requires
the coefficients of all $u^{i,j}$ and all its product terms to be equal to zero.
The determining equations for
$\xi, \tau, \eta$ are
\begin{footnotesize}
\begin{gather}
   \frac{\partial \tau  }{\partial x}=0
   \label{eq1}\\
    \frac{\partial \tau  }{\partial u}=0
    \label{eq2}\\
    \frac{\partial \xi }{\partial t}=0
    \label{eq3}\\
     \frac{\partial \xi }{\partial u}=0
     \label{eq4}\\
     m(x) \frac{\partial^{2}\eta}{\partial t^{2}}-T'(x)\frac{\partial \eta}{\partial x}-T(x)\frac{\partial^{2}\eta}{\partial x^{2}}+EI''(x)\frac{\partial^{2}\eta}{\partial x^{2}}+2EI'(x)\frac{\partial^{3}\eta}{\partial x^{3}}+EI(x)\frac{\partial^{4}\eta}{\partial x^{4}}=0
    \label{eq5}\\
    -\frac{2\xi (x,t,u)(EI'(x))^{2}}{EI(x)}+2 \xi (x,t,u)EI''(x)+2EI'(x)\frac{\partial \xi}{\partial x}+4 EI(x)\frac{\partial^{2}\eta}{\partial x \partial u}-6 EI(x)\frac{\partial^{2}\xi}{\partial x^{2}}=0
      \label{eq6}\\
    \frac{T(x)\xi (x,t,u)EI'(x)}{EI(x)}-\xi (x,t,u)T'(x)-\frac{\xi (x,t,u)EI'(x)EI''(x)}{EI(x)}+\xi (x,t,u)EI'''(x)\nonumber \\-2T(x) \frac{\partial \xi}{\partial x} +2 EI''(x)\frac{\partial \xi}{\partial x}
   +6EI'(x)\frac{\partial^{2}\eta}{\partial x \partial u}+6EI'(x)\frac{\partial^{2}\xi}{\partial x^{2}}+6EI(x)\frac{\partial^{3}\eta}{\partial x^{2}\partial u}
    +4EI(x)\frac{\partial^{3}\xi}{\partial x^{3}}=0
    \label{eq7}
\end{gather}
\begin{gather}
  \frac{T'(x)\xi (x,t,u)EI'(x)}{EI(x)}-T''(x)\xi (x,t,u)-m(x)\frac{\partial^{2}\xi}{\partial t^{2}}
  -3T'(x)\frac{\partial \xi}{\partial x}
  -2T(x)\frac{\partial^{2}\eta}{\partial x \partial u}
  +2EI''(x)\frac{\partial^{2}\eta}{\partial x \partial u}\nonumber \\+T(x)\frac{\partial^{2}\xi}{\partial x^{2}}-EI(x)\frac{\partial^{2}\xi}{\partial x^{2}}+6EI'(x)\frac{\partial^{3}\eta}{\partial x^{2}\partial u}-2EI'(x)\frac{\partial^{3}\xi}{\partial x^{3}}+4EI(x)\frac{\partial^{4}\eta}{\partial x^{3}\partial u}+EI(x)\frac{\partial^{4}\xi}{\partial x^{4}}=0\label{eq8}\\
  -\frac{m(x)\xi (x,t,u)EI'(x)}{EI(x)}+\xi (x,t,u)m'(x)-2 m(x)\frac{\partial \tau}{\partial t}+4\frac{\partial \xi}{\partial x}=0\label{eq9}\\
  2 m(x)\frac{\partial^{2}\eta}{\partial t\partial u}-m(x)\frac{\partial^{2} \tau}{\partial t^{2}}+T'(x)\frac{\partial \tau}{\partial x}+T(x)\frac{\partial^{2} \tau}{\partial x^{2}}-EI''(x)\frac{\partial^{2} \tau}{\partial x^{2}}-2EI'(x)\frac{\partial^{3} \tau}{\partial x^{3}}-EI(x)\frac{\partial^{4} \tau}{\partial x^{4}}=0\label{eq10}\\
  \frac{\partial^{2}\eta}{\partial u^{2}}=0\label{etauu}
     \end{gather}
\end{footnotesize}
From Eqs.(\ref{eq1}, \ref{eq2}) $\tau(x,t,u)=\tau(t)$, from Eqs. (\ref{eq3}, \ref{eq4}) $\xi(x,t,u)=\xi(x)$  and from Eq. (\ref{etauu})
$\eta(x,t,u)=A(x,t)u+B(x,t)$ for some arbitrary functions $A(x,t)$, $B(x,t)$. We assume $\frac{\partial B}{\partial x}=0$ and $\frac{\partial^{2}B}{\partial t^{2}}=0$, i.e., $B(x,t)=d_{1}+d_{2}t$ where $d_{1},d_{2}$ are constants. From the Eq.(\ref{eq6}), we observe that $\frac{\partial^{2}\eta}{\partial x \partial u}=\frac{\partial A}{\partial x}$ should be free from time variable $t$. Hence, $A(x,t)=f_{1}(x)+f_{2}(t)$.
\par
Now from Eq. (\ref{eq10}), we see that
\begin{equation}
2 m(x)\frac{\partial^{2}\eta}{\partial t\partial u}-m(x)\frac{\partial^{2} \tau}{\partial t^{2}}=0
\label{reduced10}
\end{equation}
Again from Eq. (\ref{eq9}), $\frac{\partial \tau}{\partial t}$ should be a constant. Assume, $2\frac{\partial \tau}{\partial t}=\omega$ which implies $\tau(t)=\frac{\omega}{2}t+t_{0}$. As from Eq. (\ref{reduced10}),
\begin{equation}
2 f_{2}(t)=\frac{\partial \tau}{\partial t}=\frac{\omega}{2}
\end{equation}
So, $f_{2}(t)=\frac{\omega}{4}$ and $\eta(x,t,u)=(f_{1}(x)+\frac{\omega}{4})u+d_{1}+d_{2}t$.
\par
Now we consider two cases;
$(a)$ when $f_{1}(x)$ is constant, $(b)$ when $f_{1}(x)$ is not constant. For the first case, there are two subcases;
$(a.1)$
when $\frac{d^{3}\xi}{dx^{3}}=k$ for some constant $k$ i.e., $\frac{d^{3}\xi}{dx^{3}}=0,\frac{d^{4}\xi}{dx^{4}}=0$,$(a.2)$ when$\frac{d^{3}\xi}{dx^{3}}$ is not a constant. All the combinations of the stiffness, mass, axial force, and transformations are given in Table \ref{symmetryTab1}.
\paragraph*{Case $(a)$ when $f_{1}(x)$ is constant.$-$}
For the case $(a.1)$, from Eq. (\ref{eq6}), it can be shown that
\begin{equation}
\xi(x)=f_{0}(EI(x))^{\frac{1}{3}}
\end{equation}
and for the assumption $\frac{d^{2}\xi}{dx^{2}}=k$, $EI(x)=\frac{k^3 x^6}{8 f_{0}^3}$.
For the case $(a.2)$, $\frac{d^{2}\xi}{dx^{2}}$ is not a constant. Here, we assume the coefficients of $\xi(x)$ are zero in the Eqs.  (\ref{eq6}), (\ref{eq7}), (\ref{eq8}) and the form of $EI(x), m(x)$, and $T(x)$ are evaluated.
\paragraph*{Case $(b)$ when $f_{1}(x)$ is not a constant.$-$}
In this case also, if we assume the coefficients of $\xi(x)$ are zero in the Eqs. (\ref{eq6}), (\ref{eq7}), (\ref{eq8}) and the form of $EI(x), m(x)$ and $T(x)$ are evaluated.
\par
Based on the assumptions, there may be more than these combinations of physical properties and coordinate transformations which can lead us to find a closed form solution. In Table \ref{symmetryTab1}, the possible list of these combinations are listed. It is observed that the transformation rule for the spatial coordinate is strongly dependent on the stiffness of the beam, i.e., the geometry of beam of certain material. It can be shown that for a beam of polynomial varying stiffness $EI(x)=(a_{0}+a_{1}x)^{n}$ where $a_{1}\neq 0$ any real number and $n> 3$ a nonzero positive integer, the proper transformation exist which leads to exact solution given in Table \ref{symmetryTab1} where $G(x)$ is
  \begin{gather}
G(x)=(a_{0}+a_{1} x)^{-n} \left(-\frac{2 \sqrt{a_{1}} A_{1} \sqrt{n-2}
   (a_{0}+a_{1} x)^{\frac{1}{2} \left(-\frac{\sqrt{a_{1}^3 (n-2)
   (n-1)^2+4  T_{1}}}{a_{1}^{3/2} \sqrt{n-2}}+n+3\right)}}{a_{1}^{3/2}
   (n-3) \sqrt{n-2}+\sqrt{a_{1}^3 (n-2) (n-1)^2+4 T_{1}}}
   -\right.\nonumber\\
    \left.
   \frac{2\sqrt{a_{1}} A_{2} \sqrt{n-2} (a_{0}+a_{1} x)^{\frac{1}{2}
   \left(\frac{\sqrt{a_{1}^3 (n-2) (n-1)^2+4  T_{1}}}{a_{1}^{3/2}
   \sqrt{n-2}}+n+3\right)}}{a_{1}^{3/2} (n-3) \sqrt{n-2}-\sqrt{a_{1}^3
   (n-2) (n-1)^2+4  T_{1}}}-
   \frac{A_{3} (a_{0}+a_{1} x)^3}{a_{1}
   (n-3)}\right)
   \label{G(x)}
\end{gather}
\begin{table}[h!]
 \centering
     \begin{tabular}{|c|l|p{3cm}|p{4cm}|}
\hline
Case & Physical Properties & Transformations & $u(x,t)$ \\ \hline
\multirow{3}{*}{$(a.1)$} & $EI(x)=\frac{k^3 x^6}{8 f_{0}^3}$ & $\xi(x)=\frac{k x^2}{2}$ &
$ e^{-2/x} \left(A_{1} e^{\lambda t}+A_{2} e^{-\lambda t}\right)
   $\\
 & $m(x)=m_{0}\frac{ e^{\frac{-4 c_{2}}{k x}}}{x^2}$ & $\tau(t)=\frac{\omega}{2}t+t_{0}$& \\
 & $T(x)=T_{0} x^2-\frac{3 k^3 x^4}{4 f_{0}^3}$ & $\eta=(\alpha+\frac{\omega}{4})u+d_{1}+d_{2}t$ &\\ \hline
\multirow{3}{*}{$(a.2)$} & $EI(x)=a_{1} e^{a_{0} x}$ & $\xi(x)=\frac{3 f_{0} e^{\frac{a_{0} x}{3}}}{a_{0}}$&
$e^{\frac{x}{r_{0}}} (A_{1}+A_{2} t) $\\
 & $m(x)=m_{0} e^{\frac{1}{3} a_{0} \left(-\frac{6 c_{2} e^{-\frac{a_{0}
   x}{3}}}{a_{0} f_{0}}-x\right)}$& $\tau(t)=\frac{\omega}{2}t+t_{0}$ &\\
 & $T(x)=-\frac{2}{9} a_{0}^2 a_{1} e^{a_{0} x}$ & $\eta=(\alpha+\frac{\omega}{4})u+d_{1}+d_{2}t$ &\\ \hline
\multirow{3}{*}{$(b)$} & $EI(x)=a_{1} e^{-v x}$ & $\xi(x)=\frac{e^{v x}}{2 v^2}$ &$e^{2 v x} (A_{1}+A_{2} t) $\\
 & $m(x)= m_{0}e^{v \left(-4 c_{2} e^{-v x}-5 x\right)}$ & $\tau(t)=\frac{\omega}{2}t+t_{0}$&\\
 & $T(x)=2 a_{1} v^2 e^{-v x}$ & $\eta=\left(\frac{e^{v x}}{v}+\frac{\omega}{4}\right)u+d_{1}+d_{2}t$ &\\ \hline
\multirow{3}{*}{$(c)$} & $EI(x)=(a_{0}+a_{1}x)^{n}$ & $\xi(x)=\frac{(a_{0}+a_{1} x)}{a_{1} n}$ &$2t+G(x)$\\
 & $m(x)=m_{0} (a_{0}+a_{1} x)^{\frac{f_{0} (n-4)+n \omega }{f_{0}}}$ & $\tau(t)=\frac{\omega}{2}t+t_{0}$ &\\
 & $T(x)=\frac{T_{1} (a_{0}+a_{1} x)^{n-2}}{a_{1} (n-2)}$ & $\eta=\frac{\omega}{4}u+d_{1}+d_{2}t$& \\
\hline
\end{tabular}
\caption{The transformations and solutions for the corresponding physical properties. Here $\lambda=-\frac{\sqrt{2} \sqrt{2 f_{0}^3
   T_{0}-k^3}}{f_{0}^{3/2} \sqrt{m_{0}}}$ and $G(x)$ is given in (\ref{G(x)}).
      }
      \label{symmetryTab1}
\end{table}
\section{Closed form solution}
For the case $(a.1)$, the acting vector field is
\begin{eqnarray}
v=\frac{k x^2}{2}\frac{\partial}{\partial x}+(\frac{\omega}{2}t+t_{0})\frac{\partial}{\partial t}+((\alpha+\frac{\omega}{4})u+d_{1}+d_{2}t)\frac{\partial}{\partial u}
\end{eqnarray}
If we choose, $(\alpha+\frac{\omega}{4})=k$, then for vector field $X=\frac{x^2}{2}\frac{\partial}{\partial x}+u\frac{\partial}{\partial u}$ the characteristic equation is
\begin{eqnarray}
\frac{dx}{\frac{x^2}{2}}=\frac{dt}{0}=\frac{du}{u}
\end{eqnarray}
which implies
\begin{eqnarray}
u(x,t)=exp\left(-\frac{2}{x}\right)F(t)
\end{eqnarray}
for some arbitrary function $F(t)$. Now substituting this $u(x,t)$ in the original equation Eq. (\ref{first equation}) for $EI(x)=\frac{k^3 x^6}{8 f_{0}^3}$, $T(x)=T_{0} x^2-\frac{3 k^3 x^4}{4 f_{0}^3}$ and $m(x)=m_{0}\frac{ e^{\frac{-4 c_{2}}{k x}}}{x^2}$, a second order ordinary differential equation in $F(t)$ is found. Solving this equation for $F(t)$, the given $u(x,t)$ in Table \ref{symmetryTab1} is found. Following the similar way, we get the invariant solutions for other cases given in the Table \ref{symmetryTab1}.
\section{Boundary condition imposition}
The symmetry analysis of a boundary value problem (BVP) may not always be successful \cite{hydon2005symmetry}. For a BVP, the domain is fixed and the symmetry of a BVP requires not only the invariance of the given differential equation, but also the invariance of the boundary data. There is no available systematic procedure for the symmetry analysis of a BVP \cite{clarkson2001open}. However, in this paper, the invariant solution is found by weakening the conditions given by Bluman and Kumei
in section $4.4.1$ \cite{bluman2013symmetries} for the invariance of the boundary.
We consider a beam with unit length. The boundary conditions are; $u(0,t)=0, \frac{\partial u}{\partial x}(0,t)=0, \frac{\partial^{2} u}{\partial x ^{2}}(1,t)=0, \frac{\partial^{3} u}{\partial x ^{3}}(1,t)=0,$ and the initial conditions are $u(x,0)=h(x)$,
$\frac{\partial u}{\partial t}(x,0)=0$.
Satisfying all the boundary condition for the deflection, i.e. for $u$, of an elastic beam which is fixed at the left end and free at the right end, the feasible stiffness, mass, and the axial load are given by;
\begin{figure}[h!]
\centering
\includegraphics[scale=0.7]{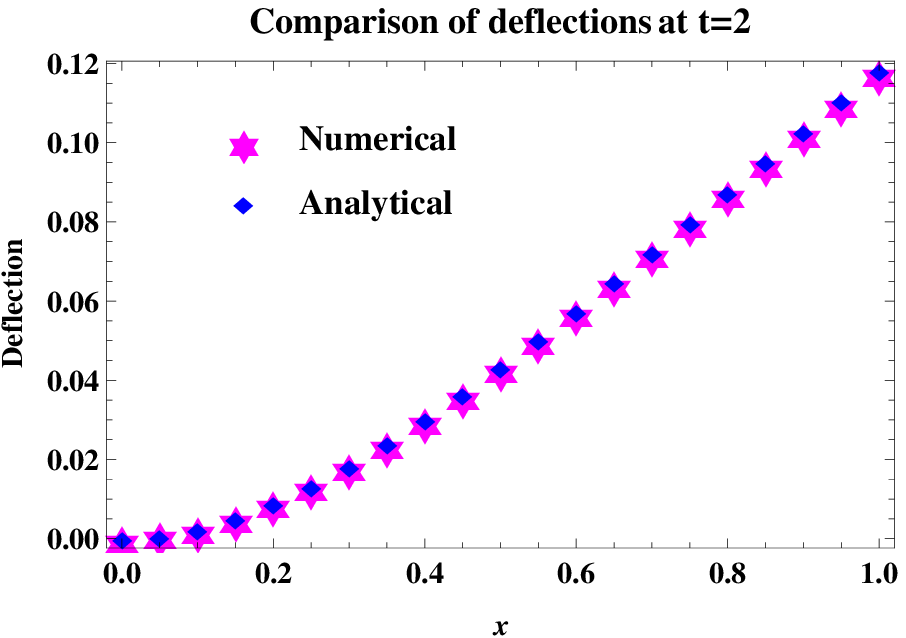}
\includegraphics[scale=0.7]{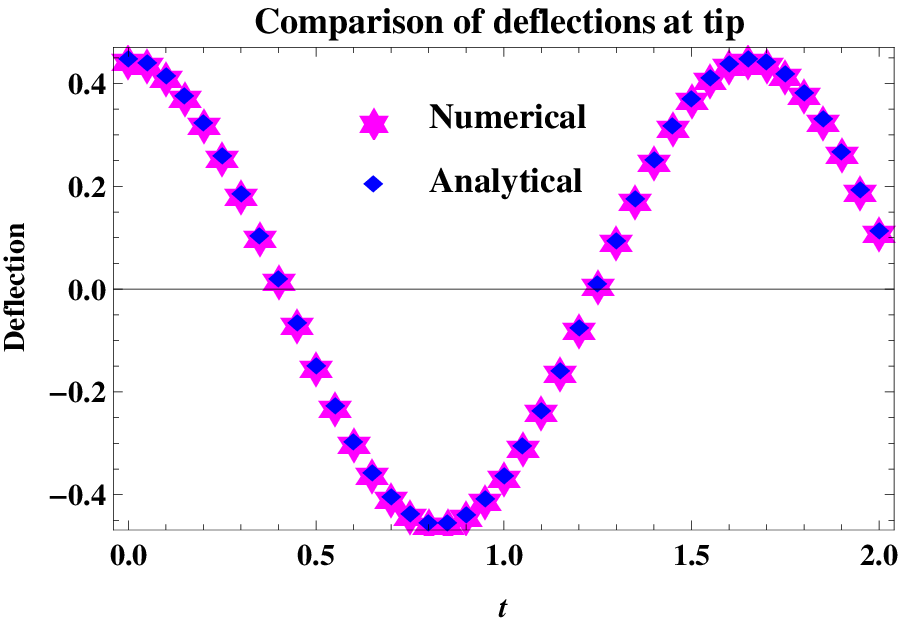}
\caption{Comparison of numerical and analytical solutions }
\label{comparison}
\end{figure}
\begin{eqnarray}
EI(x)=\frac{g_{1} x^4 (6 m_{0}+x (3 m_{1}+2 m_{2} x))^4}{(m_{0}+x
   (m_{1}+m_{2} x))^3}
   \end{eqnarray}
   \begin{eqnarray}
 m(x)=  m_{0}+m_{1} x+m_{2} x^2
 \end{eqnarray}
 \begin{gather}
 T(x)=\frac{g_{1} x^2 (6 m_{0}+x (3 m_{1}+2 m_{2} x))^2 }{(m_{0}+x (m_{1}+m_{2} x))^5}
 \left(-56
   m_{0}^4+8 m_{0}^3 x (17 m_{2} x-10 m_{1})
    \right.\nonumber\\\left.
+12 m_{0}^2 x^2
   \left(-7 m_{1}^2+8 m_{1} m_{2} x+2 m_{2}^2 x^2\right)
   +2m_{0} x^3 \left(-22 m_{1}^3+9 m_{1}^2 m_{2} x+24 m_{1}
   m_{2}^2 x^2+12 m_{2}^3 x^3\right)
   \right.\nonumber\\\left.
   -m_{1}^2 x^4 \left(11
   m_{1}^2+14 m_{1} m_{2} x+6 m_{2}^2
   x^2\right)\right)
\end{gather}
and the transformations are;
\begin{eqnarray}
\xi(x)=\frac{m_{0} x+\frac{m_{1} x^2}{2}+\frac{m_{2} x^3}{3}}{m_{0}+x(m_{1}+m_{2} x)}\\
   \tau(t)=c_{1}\\
   \eta=\omega u+d_{1}+d_{2}t
   \end{eqnarray}
where $m_{1}, m_{2}$ are dependent on $m_{0}$ and the general expressions are
\begin{equation*}
m_{1}=\frac{\sqrt[3]{\sqrt{6} m_{0}^3+9m_{0}^3}}{3^{2/3}
   \sqrt[3]{5}}+\frac{\sqrt[3]{\frac{5}{3}} m_{0}^2}{\sqrt[3]{\sqrt{6}
   m_{0}^3+9 m_{0}^3}}-2 m_{0}
\end{equation*}
$
m_{2}=  \frac{1}{240 m_{0}^2}( 5 \sqrt{6} m_{0}^3+195 m_{0}^3+8 \sqrt[3]{5} \left(3
   \left(9+\sqrt{6}\right)\right)^{2/3} \left(m_{0}^3\right)^{2/3}
   m_{0}+\frac{375 m_{0}^6}{\sqrt{6} m_{0}^3+9 m_{0}^3}+\frac{120
   \sqrt[3]{3} 5^{2/3} m_{0}^5}{\left(9+\sqrt{6}\right)^{2/3}
   \left(m_{0}^3\right)^{2/3}}-\frac{120\ 3^{2/3}
   \sqrt[3]{\frac{5}{9+\sqrt{6}}} m_{0}^4}{\sqrt[3]{m_{0}^3}}-24\ 5^{2/3}
   \sqrt[3]{3 \left(9+\sqrt{6}\right)} \sqrt[3]{m_{0}^3} m_{0}^2)
$
\par
We choose $\omega = 2$, $g_{0} = 1$, $g_{1} = \frac{1}{3000}$, $m_{0}=1$. Then $m_{1}=-0.840295, m_{2}=0.277816$. Here $\xi(0)=0$ but it is not possible to impose $\xi(1)=0$ which provides non-feasible physical properties.
The corresponding analytical solution is ;
 \begin{eqnarray}
u(x,t)=\left(m_{0} x+\frac{m_{1} x^2}{2}+\frac{m_{2} x^3}{3}\right)^2
   \left(A_{1} \cos \left(120 \sqrt{3} \sqrt{g_{1}} t\right)+A_{2} \sin
   \left(120 \sqrt{3} \sqrt{g_{1}} t\right)\right)
    \end{eqnarray}
For zero initial velocity, $A_{2} =0$. Considering $A_{1}=1$,
 \begin{eqnarray}
u(x,t)=\left(m_{0} x+\frac{m_{1} x^2}{2}+\frac{m_{2} x^3}{3}\right)^2
   \left( \cos \left(120 \sqrt{3} t\right)+ \sin
   \left(120 \sqrt{3}  t\right)\right)
    \end{eqnarray}
    The comparison between the analytical form of the solution is matching well with the numerical solution given in Fig. \ref{comparison}. For the numerical solution, NDSolve command in $Mathematica$ $9.0$ is used with $h(x)=\left(m_{0} x+\frac{m_{1} x^2}{2}+\frac{m_{2} x^3}{3}\right)^2$.
\section{Conclusion}
In conclusion, Lie symmetry method is applied to analyze the symmetry and the exact or closed form solution of the Euler-Bernoulli beam with axial load. Some combinations of the coordinate transformations dependent on the system properties are found which provide the closed form solution.
Different combinations of stiffness, mass and axial force are also found which yield
a closed form solution. It is observed that the crucial spatial transformations are dependent on the stiffness of the beam.
There are some nanobeam type structures such as nanoropes, nanorod etc. \cite{fan2009superelastic,lu1997elastic,manghi2006propulsion} where this analytical result can be used directly to measure the strain-displacement relation.
The imposition of boundary conditions for the symmetry analysis of an elastic beam is carried out successfully.
\bibliographystyle{ieeetr}
\bibliography{thesis_bib}

\begin{thebibliography}{10}

\bibitem{guede2001apparentlynew}
Z.~Guede and I.~Elishakoff, ``Apparently first closed-form solutions for
  inhomogeneous vibrating beams under axial loading,'' {\em Proc. R. Soc. A},
  vol.~457, no.~2007, pp.~623--649, 2001.

\bibitem{ovsjannikov1962group}
L.~Ovsjannikov and G.~Bluman, {\em Group properties of differential equations}.
\newblock Siberian Section of the Academy of Science of USSR, 1962.

\bibitem{muller1962and}
E.~M{\"u}ller and K.~Matschat, ``{\"U}ber das auffinden von
  {\"a}hnlichkeitsl{\"o}sungen partieller differentialgleichungssysteme unter
  ben{\"u}tzung von transformationsgruppen,'' 1962, p 190.

\bibitem{olver2012applications}
P.~J. Olver, {\em Applications of Lie Groups to Differential Equations},
  vol.~107.
\newblock Springer Science \& Business Media, New York, 2012.

\bibitem{bluman2013symmetries}
G.~Bluman and S.~Kumei, {\em Symmetries and differential equations}, vol.~154.
\newblock Springer Science \& Business Media, New York, 2013.

\bibitem{ibragimov1995crc}
N.~H. Ibragimov, {\em CRC Handbook of Lie group analysis of differential
  equations}, vol.~3.
\newblock CRC press,Boca Raton, FL, 1995.

\bibitem{hydon2000symmetrynew}
P.~E. Hydon, {\em Symmetry methods for differential equations: a beginner's
  guide}, vol.~22.
\newblock Cambridge University, New York, 2000.

\bibitem{bluman1980remarkable}
G.~Bluman and S.~Kumei, ``On the remarkable nonlinear diffusion equation,''
  {\em Journal of Mathematical Physics}, vol.~21, no.~5, pp.~1019--1023, 1980.

\bibitem{bluman1987invariance}
G.~Bluman and S.~Kumei, ``On invariance properties of the wave equation,'' {\em
  Journal of mathematical physics}, vol.~28, no.~2, pp.~307--318, 1987.

\bibitem{torrisi1996group}
M.~Torrisi, R.~Tracina, and A.~Valenti, ``A group analysis approach for a
  nonlinear differential system arising in diffusion phenomena,'' {\em Journal
  of Mathematical Physics}, vol.~37, no.~9, pp.~4758--4767, 1996.

\bibitem{ibragimov2004equivalence}
N.~H. Ibragimov and N.~S{\"a}fstr{\"o}m, ``The equivalence group and invariant
  solutions of a tumour growth model,'' {\em Communications in Nonlinear
  Science and Numerical Simulation}, vol.~9, no.~1, pp.~61--68, 2004.

\bibitem{ibragimov2011lie}
N.~H. Ibragimov, ``Lie group analysis of {M}offatt's model in metallurgical
  industry,'' {\em Journal of Nonlinear Mathematical Physics}, vol.~18,
  no.~sup1, pp.~143--162, 2011.

\bibitem{gray2015calculate}
R.~J. Gray, ``How to calculate all point symmetries of linear and linearizable
  differential equations,'' {\em Proc. R. Soc. A}, vol.~471, no.~2175,
  p.~20140685, 2015.

\bibitem{gainetdinova2017integrability}
A.~Gainetdinova and R.~Gazizov, ``Integrability of systems of two second-order
  ordinary differential equations admitting four-dimensional lie algebras,''
  {\em Proc. R. Soc. A}, vol.~473, no.~2197, p.~20160461, 2017.

\bibitem{hydon2000symmetries}
P.~Hydon, ``Symmetries and first integrals of ordinary difference equations,''
  {\em Proc. R. Soc. A}, vol.~456, no.~2004, pp.~2835--2855, 2000.

\bibitem{kang2012symmetry}
J.~Kang and C.~Qu, ``Symmetry groups and fundamental solutions for systems of
  parabolic equations,'' {\em Journal of Mathematical Physics}, vol.~53, no.~2,
  p.~023509, 2012.

\bibitem{singla2017invariant}
K.~Singla and R.~Gupta, ``On invariant analysis of space-time fractional
  nonlinear systems of partial differential equations. ii,'' {\em Journal of
  Mathematical Physics}, vol.~58, no.~5, p.~051503, 2017.

\bibitem{hau2017optimal}
J.-N. Hau, M.~Oberlack, and G.~Chagelishvili, ``On the optimal systems of
  subalgebras for the equations of hydrodynamic stability analysis of smooth
  shear flows and their group-invariant solutions,'' {\em Journal of
  Mathematical Physics}, vol.~58, no.~4, p.~043101, 2017.

\bibitem{strobl2004algebroid}
T.~Strobl, ``Algebroid yang-mills theories,'' {\em Physical Review Letters},
  vol.~93, no.~21, p.~211601, 2004.

\bibitem{beisert2017yangian}
N.~Beisert, A.~Garus, and M.~Rosso, ``Yangian symmetry and integrability of
  planar n= 4 supersymmetric yang-mills theory,'' {\em Physical Review
  Letters}, vol.~118, no.~14, p.~141603, 2017.

\bibitem{belmonte2007lie}
J.~Belmonte-Beitia, V.~M. P{\'e}rez-Garc{\'\i}a, V.~Vekslerchik, and P.~J.
  Torres, ``Lie symmetries and solitons in nonlinear systems with spatially
  inhomogeneous nonlinearities,'' {\em Physical Review Letters}, vol.~98,
  no.~6, p.~064102, 2007.

\bibitem{budanur2015reduction}
N.~B. Budanur, P.~Cvitanovi{\'c}, R.~L. Davidchack, and E.~Siminos, ``Reduction
  of so (2) symmetry for spatially extended dynamical systems,'' {\em Physical
  Review Letters}, vol.~114, no.~8, p.~084102, 2015.

\bibitem{bocko2012symmetries}
J.~Bocko, V.~Nohajov{\'a}, and T.~Har{\v{c}}arik, ``Symmetries of differential
  equations describing beams and plates on elastic foundations,'' {\em Procedia
  Engineering}, vol.~48, pp.~40--45, 2012.

\bibitem{ozkaya2002group}
E.~{\"O}zkaya and M.~Pakdemirli, ``Group--theoretic approach to axially
  accelerating beam problem,'' {\em Acta Mechanica}, vol.~155, no.~1-2,
  pp.~111--123, 2002.

\bibitem{soh2008euler}
C.~W. Soh, ``Euler-{B}ernoulli beams from a symmetry
  standpoint-characterization of equivalent equations,'' {\em Journal of
  Mathematical Analysis and Applications}, vol.~345, no.~1, pp.~387--395, 2008.

\bibitem{bokhari2010symmetries}
A.~H. Bokhari, F.~Mahomed, and F.~Zaman, ``Symmetries and integrability of a
  fourth-order {E}uler-{B}ernoulli beam equation,'' {\em Journal of
  Mathematical Physics}, vol.~51, no.~5, p.~053517, 2010.

\bibitem{johnpillai2016noether}
A.~Johnpillai, K.~Mahomed, C.~Harley, and F.~Mahomed, ``Noether symmetry
  analysis of the dynamic {E}uler-{B}ernoulli beam equation,'' {\em Zeitschrift
  f{\"u}r Naturforschung A}, vol.~71, no.~5, pp.~447--456, 2016.

\bibitem{gunda2008stiff}
J.~B. Gunda and R.~Ganguli, ``Stiff-string basis functions for vibration
  analysis of high speed rotating beams,'' {\em Journal of Applied Mechanics},
  vol.~75, no.~2, pp.~024502--024506, 2008.

\bibitem{hassani2013mathematical}
S.~Hassani, {\em Mathematical physics: a modern introduction to its
  foundations}.
\newblock Springer Science \& Business Media, 2013.

\bibitem{hydon2005symmetry}
P.~E. Hydon, ``Symmetry analysis of initial-value problems,'' {\em Journal of
  Mathematical Analysis and Applications}, vol.~309, no.~1, pp.~103--116, 2005.

\bibitem{clarkson2001open}
P.~A. Clarkson and E.~L. Mansfield, ``Open problems in symmetry analysis,''
  {\em Contemporary Mathematics}, vol.~285, pp.~195--205, 2001.

\bibitem{fan2009superelastic}
W.~Fan, S.~Huang, J.~Cao, E.~Ertekin, C.~Barrett, D.~Khanal, J.~Grossman, and
  J.~Wu, ``Superelastic metal-insulator phase transition in single-crystal vo 2
  nanobeams,'' {\em Physical Review B}, vol.~80, no.~24, p.~241105, 2009.

\bibitem{lu1997elastic}
J.~P. Lu, ``Elastic properties of carbon nanotubes and nanoropes,'' {\em
  Physical Review Letters}, vol.~79, no.~7, p.~1297, 1997.

\bibitem{manghi2006propulsion}
M.~Manghi, X.~Schlagberger, and R.~R. Netz, ``Propulsion with a rotating
  elastic nanorod,'' {\em Physical Review Letters}, vol.~96, no.~6, p.~068101,
  2006.

\end{thebibliography}
\end{document}